\documentclass[11pt]{article}

\usepackage[utf8]{inputenc}
\usepackage{latexsym}
\usepackage[T1]{fontenc}
\usepackage[english]{babel}

\usepackage[a4paper, margin=2.5cm]{geometry}
\usepackage{setspace}
\usepackage{amsmath, amssymb, amsthm}
\usepackage{amsfonts}

\usepackage{graphicx}
\usepackage{caption}
\usepackage{subcaption}
\usepackage{booktabs}
\usepackage{float}

\usepackage{xcolor}
\usepackage{hyperref}
\hypersetup{
    colorlinks=true,
    linkcolor=blue,
    citecolor=blue,
    urlcolor=blue
}

\title{Designing for the Development of Probabilistic Thinking: A Design-Based Research Study in Lower Secondary Education}
\author{Luigia Caputo$^{1,*}$, Aniello Buonocore$^2$}
\date{\small Università di Napoli Federico II -- Dipartimento di Matematica e Applicazione ``Renato Caccioppoli''\\
Via Cintia, Complesso di Monte S. Angelo, 80126 Napoli, Italia\\
$^1$luigia.caputo@unina.it, $^2$aniello.buonocore@unina.it, $^*$corresponding author}

\theoremstyle{definition}

\newtheorem{remark}{Remark}

\begin{document}

\maketitle

\begin{abstract}
Drawing on the Data and Predictions strand of the Indicazioni Nazionali per il curricolo 2012, this study proposes a problem-based instructional approach to the teaching of probability.
More specifically, the study adopts a design-based research methodology structured in a single cycle consisting of two teaching interventions in the same class, carried out in two consecutive years. Within this framework, a set of carefully selected problems 
is employed to foster students’ engagement. These problems are designed not only to introduce probabilistic concepts, but also to stimulate students’ communicative and argumentative skills.
The selected tasks provide opportunities to promote key \emph{process goals} (such as reasoning and proving, communicating, representing, and making connections) which are often overshadowed by a predominant focus on \emph{content goals}. This approach aims to support teachers in addressing the difficulties they frequently encounter in guiding students’ conceptualization processes, particularly in bridging the gap between students’ intuitive reasoning and formal abstraction. At the same time, it seeks to help students develop more robust and flexible forms of thinking, enabling them to better navigate situations involving uncertainty in everyday life.
\end{abstract}
%
%
\section{Introduction}
Design research has emerged as a prominent methodological approach in educational studies, particularly suited to addressing complex problems situated in real-world learning environments. Unlike traditional research paradigms that often prioritize either theoretical abstraction or empirical description, design research seeks to establish a productive interplay between theory and practice (\cite{Bro92}, \cite{Col92}). It is increasingly recognized as a powerful approach for bridging the gap between educational theory and classroom practice, as highlighted by Design Research in Education (\cite{Bak18}).
\\
A defining characteristic of design research is its iterative nature, unfolding through cycles of design, implementation, analysis, and refinement. This cyclic process allows both educational interventions and theoretical insights to evolve progressively over time (\cite{DBRC03}, \cite{Bak18}). Such an approach is particularly relevant in domains like probability education, where learning processes are complex and often resistant to linear instructional strategies.
\\
Design research is closely connected to qualitative research traditions, which provide important methodological and epistemological underpinnings. According to Merriam and Tisdell (\cite{Mer16}), qualitative research aims to understand how individuals construct meaning from their experiences with in specific social and contextual settings. It is characterized by an interpretive stance, a naturalistic orientation, and the use of rich, descriptive data sources such as observations, interviews, and document analysis. Rather than focusing on statistical generalization, qualitative research seeks to develop in-depth understandings of phenomena from the participants' perspectives. In educational contexts, this approach is particularly suited to capturing the complexity of classroom dynamics and learning processes as they occur in authentic settings, thereby offering a foundation that aligns closely with the goals and methods of design-based research.
\\
The present study draws on both traditions and adopts a design-oriented approach informed by qualitative research principles. However, it does not fully follow the canonical iterative structure of design-based research and does not aim to develop a second cycle of refinement. Instead, the study proposes a methodological approach to elicit students' spontaneous probabilistic thinking. The goal is to capture these initial forms of reasoning and support their early formalization. In this sense, the study lies at the intersection of design and interpretation, using instructional interventions as a means to access and analyse emerging probabilistic thinking.
\\
In contemporary society, probabilistic competencies have become essential for navigating a world increasingly shaped by data, uncertainty, and random phenomena. Citizens are required to move beyond deterministic thinking and develop reasoning strategies that support informed decision-making in everyday and professional contexts. Despite the inclusion of probability in school and university curricula, research consistently shows that both students and adults — including those with formal instruction — exhibit poor probabilistic reasoning, often characterized by misconceptions and cognitive biases.
\\
In the Italian context, these issues are explicitly addressed in the \lq\lq Indicazioni Nazionali per il curricolo 2012\rq\rq\ (\cite{IN12}), where the strand data and predictions (dati e previsioni) is recognized as a key component of mathematics education from early childhood and primary education onwards. The guidelines emphasize the importance of developing students’ ability to interpret data, reason under uncertainty, and make informed predictions. However, in classroom practice, this domain is often perceived as problematic. One possible reason is that teachers may feel insufficiently confident with statistical and probabilistic content, which can lead to a marginalization of these topics or to their treatment in a purely procedural way (\cite{Bat04}).
\\
From a different perspective, however, this area can be seen as a valuable educational resource. As argued by Giorgio Arrigo (\cite{Arr10}), engaging students in probabilistic situations offers opportunities to address not only conceptual and algorithmic aspects of learning, but also other equally important dimensions. In particular, probabilistic activities can foster strategic competence, through the exploration of genuine problems in a-didactical situations; communicative competence, by encouraging students to argue, justify, and critically evaluate different viewpoints; and semiotic competence, by involving multiple representations and the ability to move between different semiotic registers.
\\
These considerations further highlight the relevance of design research for probability education. By supporting the development and testing of educational artifacts—such as tasks, simulations, and instructional sequences — in authentic classroom settings, it enables researchers to explore how students develop probabilistic understanding in practice (\cite{Ede02}, \cite{Kel03}). At the same time, it facilitates the generation of theoretical insights into learning processes.
\\
Moreover, research in probability education has expanded significantly in recent decades, as evidenced by contributions presented at major international conferences such as ICME, CERME, PME, and ICOTS, and by a growing body of literature synthesizing theoretical and empirical findings. This body of work highlights the need for innovative teaching approaches that address students’ misconceptions, support conceptual understanding, and integrate tools such as simulations, visualizations, and exploratory activities.
\\
In this context, design research offers a flexible yet rigorous framework for tackling complex and ill-structured educational problems. It promotes collaboration among researchers, teachers, and other stakeholders, and emphasizes the importance of context-sensitive solutions rather than universal prescriptions (\cite{Bar04}, \cite{Bak18}). By balancing methodological rigor with adaptability, design research enables the systematic study of teaching and learning processes while remaining responsive to the realities of classroom practice.
%
%
%
%
%
\section{Research questions and aims}
The present study aims to investigate the design and implementation of instructional activities to support the development of probabilistic thinking in lower secondary education through a design-based research approach. To address this aim, the following research questions are formulated.
\begin{enumerate}
\item What are the key characteristics of mathematical practices that can effectively support the introduction of the concept of probability in lower secondary education?
\item What design principles and pedagogical considerations should guide the development of mathematical activities that foster the formalization of probabilistic thinking among lower secondary school students?
\end{enumerate}
More specifically, the study pursues the following objectives:
\begin{itemize}
    \item[-] To design, implement, and iteratively refine a sequence of probabilistic learning activities, conducted over two consecutive school years within the same class, aimed at introducing and developing students’ understanding of probability.   
    \item[-] To explore how students engage with and make sense of informal probabilistic concepts, such as ``more likely'' and ``less likely,'' within game-based contexts during the initial phase of instruction.
    \item[-] To identify key characteristics of mathematical practices and derive design principles for tasks and learning environments that support the introduction and gradual formalization of students’ probabilistic reasoning.
\end{itemize}
%
%
%
%
%
\section{Methods and implementation}
The present study originates from an experience developed within a teacher training course organized by the Naples section of the Fondazione I Lincei per la Scuola, addressed to lower secondary school mathematics teachers.
\par
Within this context, a series of laboratory-based activities was designed and proposed with the aim of engaging students in doing mathematics, with particular reference to probability. The design of these activities was informed by prior reflections on the strand data and predictions (dati e previsioni) as outlined in the Indicazioni Nazionali per il curricolo 2012, which emphasizes the importance of developing students’ ability to reason with data and uncertainty (see, for examples, \cite{Bat16} and \cite{Bor96}).
\par
The activities were centered on the use of concrete artifacts, intended to allow students to explore probabilistic situations through hands-on experimentation. By directly interacting with these materials, students were encouraged to actively construct meaning and develop their understanding of probabilistic concepts.
\par
A fundamental assumption underpinning this work is that engagement with probabilistic situations can act as a gateway to broader mathematical practices. In particular, probability offers a productive context that can serve as a connecting bridge across different mathematical and scientific domains, supporting the development of transversal competencies and fostering meaningful links between concepts.
\\
During a lesson, a teacher from a lower secondary school in Pozzuoli Pozzuoli (near Naples, southern Italy) showed particular interest in the proposed activities, to the extent of inviting us to conduct a session with her first-year students. This opportunity made it possible to bring the designed activities into an authentic classroom context, marking the starting point for the study presented in this paper.
\par
The design of the laboratory was inspired by the Thinking Classroom approach introduced by Peter Liljedahl (\cite{Lil20}). This approach aims to create classroom environments in which students are actively engaged in thinking and problem solving, rather than passively receiving information. It emphasizes the use of collaborative work, rich and non-routine tasks, vertical non-permanent workspaces, and strategies that promote student autonomy and discussion. Within this framework, learning is seen as emerging from students’ active participation in mathematical practices, making it particularly suitable for the exploration of probabilistic situations.
\par
In particular, with respect to the selection of “good problems” to be proposed to students, initial attention was devoted to identifying tasks that were both engaging and cognitively demanding. To this end, we drew on adapted — sometimes simplified — versions of historical problems in probability, such as the well-known game of Zara (see Remark~\ref{Remark_Zara}). These tasks were chosen for their potential to stimulate curiosity and to introduce probabilistic reasoning through meaningful and accessible contexts.
\\
With regard to the learning environment, instead of having students remain seated at their desks, they were invited to stand and work on large A3 sheets attached to the walls using masking tape. This choice was intended to promote active participation, visibility of students’ work, and collective discussion, in line with the principles of the Thinking Classroom approach proposed by Peter Liljedahl.
\\
Finally, working groups were formed in a visibly random way, using playing cards to assign students to groups. This strategy was adopted to encourage collaboration, reduce fixed group dynamics, and foster a classroom culture oriented toward shared problem solving and mathematical discussion.
\begin{remark}\label{Remark_Zara}
The Game of Zara is played by throwing three dice simultaneously, and the relevant quantity is the sum of the scores obtained on each die. It has ancient origins; it is also mentioned by Dante in the first 12 lines of the sixth canto of Purgatorio in the Divine Comedy (around 1316).
\begin{quote}
Le anime, alla stregua dei postulanti che fanno cerchio intorno a chi esce vittorioso dal gioco della zara, si affollano intorno a Virgilio e a Dante, invocando da quest'ultimo preghiere e suffragi perché venga abbreviata la loro permanenza nell'Anti\-purgatorio.\footnote{The souls, much like the players who gather around the winner of the Game of Zara, crowd around Virgil and Dante, imploring the latter for prayers and suffrages to shorten their stay in Ante-Purgatory.}
\end{quote}
A few centuries later, between 1612 and 1623, an unknown patron (traditionally identified as Cosimo II de’ Medici), having observed through his own gaming experience and that of other Florentine noblemen that, despite
\begin{eqnarray*}
11&=1+4+6=1+5+5=2+3+6=2+4+5=3+3+5=3+4+4,\\
12&=1+5+6=2+4+6=2+5+5=3+3+6=3+4+5=4+4+4,
\end{eqnarray*}
the sums 11 and 12 did not occur with the same frequency, turned to Galileo Galilei to seek an explanation. Galileo published his solution in the treatise Sopra le scoperte dei dadi. Here one reads:
\begin{quote}
Ora io per servire a chi mi ha comandato che io deva produr ciò che sopra tal difficoltà mi sovviene, esporrò il mio pensiero con speranza non solamente di scior questo dubbio ma di aprir la strada a poter puntualissimamente scorger le ragioni per le quali tutte le particolarità del giuoco sono state con grande avvedimento e giudizio comparite e aggiustate.\footnote{Now, in order to serve the one who has requested that I set forth what occurs to me regarding this difficulty, I shall present my reasoning with the hope not only of resolving this doubt, but also of opening the way to clearly discerning the reasons why all the particular features of the game have been arranged and adjusted with great foresight and judgment.}
\end{quote}
Figure~\ref{Galileo} presents the solution method proposed by Galileo Galilei in his treatise.
\begin{figure}
    \centering
\includegraphics[width=12.0cm]{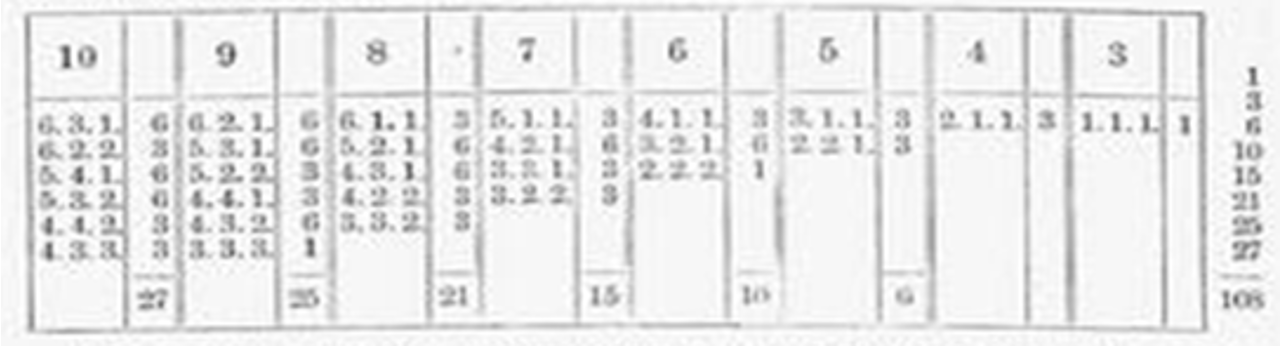}
\caption{Solution scheme proposed by Galileo Galilei in \lq\lq Sopra le scoperte dei dadi\rq\rq.}
\label{Galileo}
\end{figure}
\end{remark}
%
%
%
%
\section{First year laboratory design}
The laboratory was conducted in a first-year class of lower secondary school at the Istituto Comprensivo Giacinto Diano in Pozzuoli (Naples) in February 2020 and lasted two hours. The external instructors were present in the classroom and coordinated the laboratory activities, while the regular mathematics teacher also attended and provided logistical support. The students had not previously encountered topics related to probability. The group work was structured as cooperative learning, and the main aim of the laboratory was to allow the concept of probability to emerge naturally from the activities and to support its initial formalization.
\begin{itemize}
\item[-] The class was divided into groups of four students in a visibly random manner, labeled 1, 2, 3, and 4. 
\item[-] Each group was provided with a specific tool (a die; a set of playing cards from the same suit; or two sets of playing cards from the same suit) and a poster sheet with a marker.
\item[-] Each group was further divided into two subgroups: subgroup A (students 1 and 2) and subgroup B (students 3 and 4).
\item[-] The problem under investigation was presented and discussed with the whole class.
\item[-] The two students in subgroup A worked together for a few minutes, supported by student 3, who was responsible for handling the tool.
\item[-] One student from subgroup A was allowed to briefly consult the posters of other groups.
\item[-] The students in subgroup A then explained their solution strategy to the students in subgroup B.
\item[-] Student 3 from subgroup B implemented and simultaneously tested the proposed solution procedure.
\item[-] Student 4 from subgroup B summarized the procedure on the poster.
\end{itemize}
%
%
%
%
%
\subsection{The \lq\lq timed\rq\rq\ horse race}
The activity designed and proposed within this study is described below.
\par
\begin{figure}
    \centering
\includegraphics[width=3.5cm]{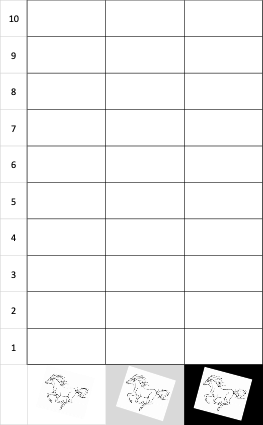}
\caption{The game board with a drawing of the three horses.}
\label{Tabellone_corsa}
\end{figure}

{\bf Activity 1}
\\
The task is based on a game scenario involving a race among three horses. The race director is provided with a box containing two dice, one blue and one red. At each round, the box is shaken and the outcomes of the two dice are observed.
The movement of the horses depends on the parity of the outcomes: the white horse advances by one square when both outcomes are even, the black horse advances when both outcomes are odd, and the grey horse advances when one outcome is even and the other is odd. The game is repeated for ten rounds.
\\
{\bf Question}: After 10 rounds, which horse has more possibilities to be ahead of the other two? Explain your reasoning.
\\
{\bf Materials}: Each group is provided with A3 sheets, markers, a game board (see Fig.~\ref{Tabellone_corsa}), and tokens representing the three horses, which support both the execution of the activity and the representation of students’ strategies and results.
\par
{\bf Activity 2}
\\
{\bf Question}: The race is repeated several times (each race consisting of ten dice throws), and the results are recorded in a table. Students are then invited to reflect on the outcomes obtained and to reconsider the identification of the horse in the previous activity.
\\
{\bf Materials}: Each group is provided with A3 sheets, markers, a cardboard game board, a box containing the two dice, and tokens representing the three horses.
\par
{\bf Activity 3}
\\
{\bf Question}: How can the rules of the game be modified so that the black and white horses have the same probability of winning as the grey horse?
%
%
%
%
%
\subsection{Reflections on the results of the first laboratory}
\begin{figure}
    \centering
\includegraphics[width=7.5cm]{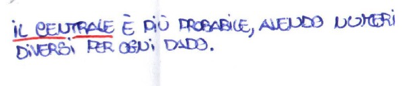}
\hspace{1cm}
\includegraphics[width=3.5cm]{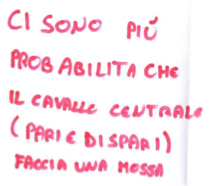}
\caption{On the left, a written remark by the group labeled A: \lq\lq The middle one is more likely, since each die shows a different type of number\rq\rq. On the right, a written remark by the group labeled B: \lq\lq There is a higher probability that the middle horse (even and odd) will make a move\rq\rq.}
\label{Cavalli_gruppo1_2}
\end{figure}

\begin{figure}
    \centering
\includegraphics[width=3.5cm]{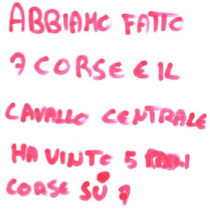}
\caption{\lq\lq We carried out 7 races, and the middle horse won 5 out of 7 races\rq\rq.}
\label{Frequenza_Probabilità}
\end{figure}

\begin{figure}
    \centering
\includegraphics[width=5.5cm]{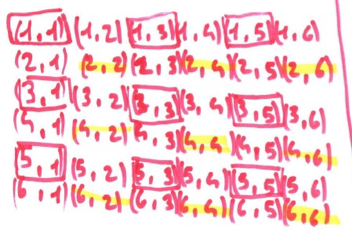}
\hspace{1cm}
\includegraphics[width=5.5cm]{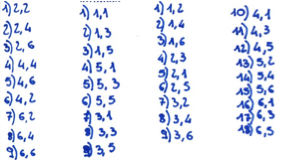}
\caption{Groups’ organization of the pairs of possible outcomes.}
\label{Coppie_gruppo1_2}
\end{figure}

\begin{figure}
    \centering
\includegraphics[width=3.5cm]{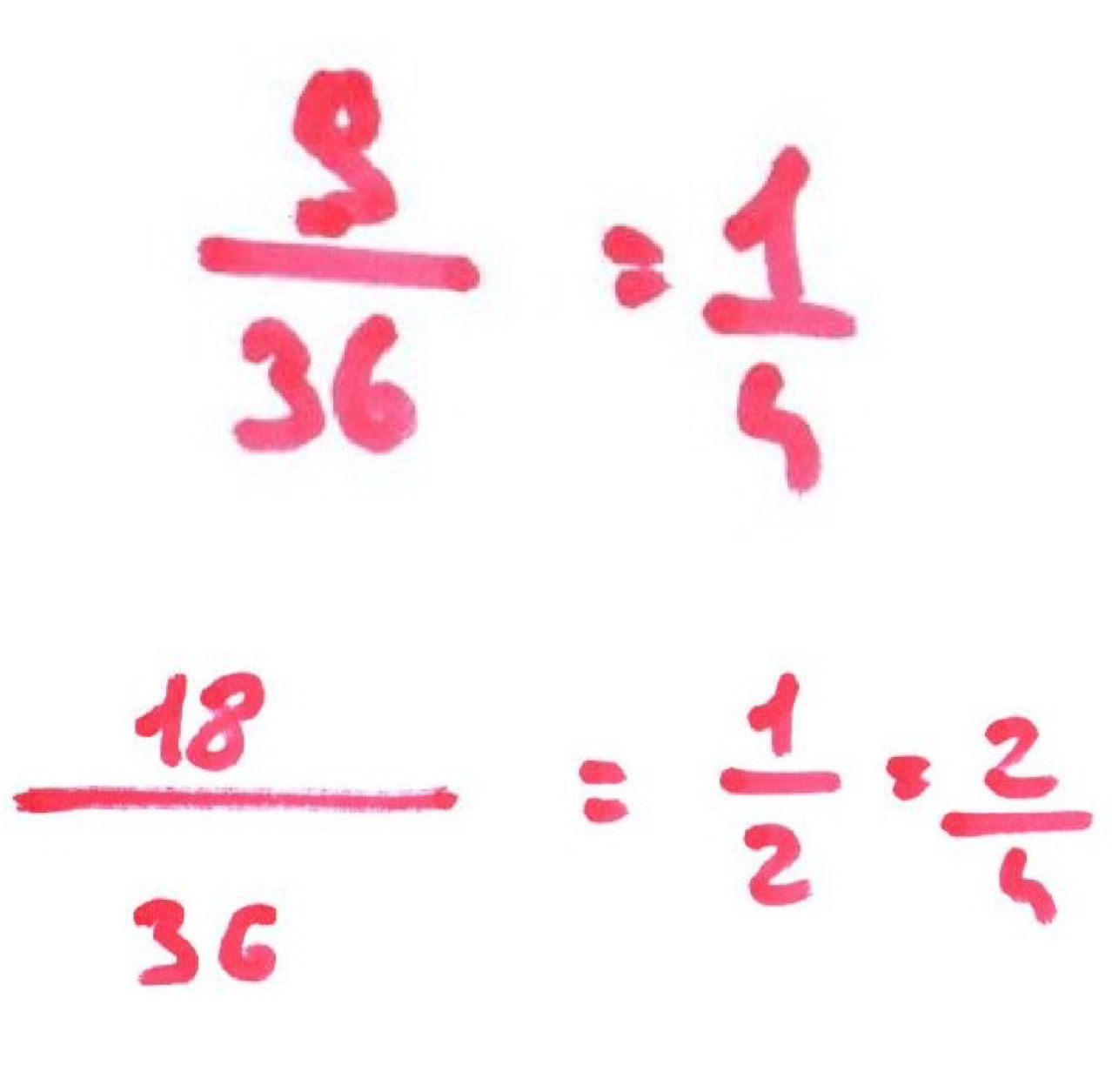}
\hspace{1cm}
\includegraphics[width=8.5cm]{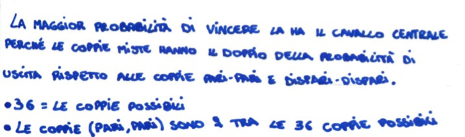}
\caption{Probability calculations carried out by some student groups.}
\label{Valutazione_probabilità}
\end{figure}

\begin{figure}
    \centering
\includegraphics[width=7.0cm]{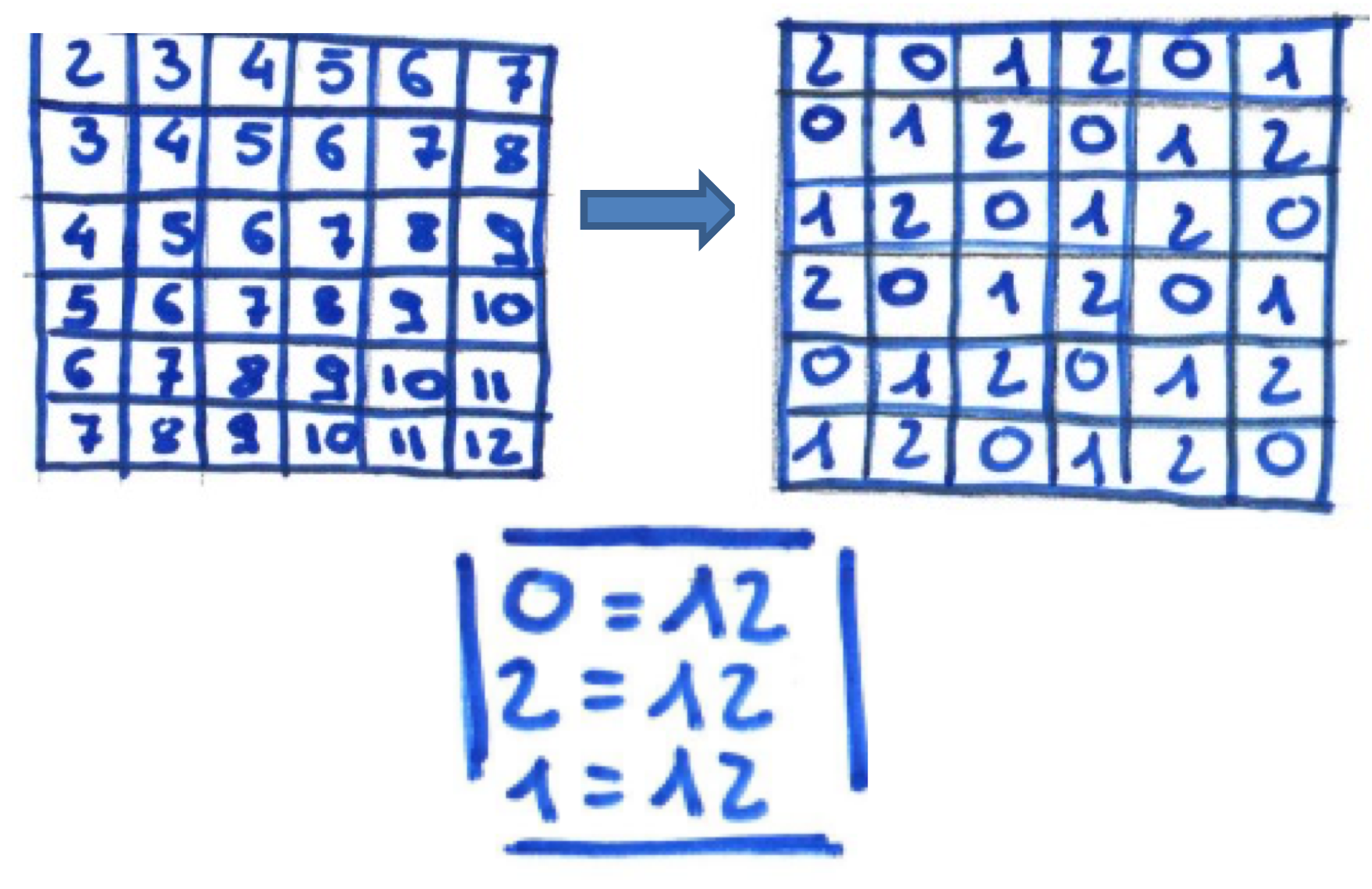}
\hspace{1cm}
\includegraphics[width=3.5cm]{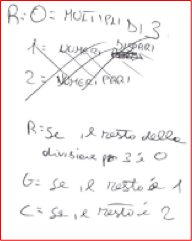}
\caption{Reorganizing data according to the modulo 3 operation.}
\label{modulo3}
\end{figure}

With reference to Activity 1, in each group there was at least one student who immediately oriented their reasoning in the correct direction and attempted to involve the other group members. In particular, these students recognized that the pairs consisting of one even and one odd outcome are more numerous than those with outcomes of the same parity.
\\
More challenging for the students was the systematic organization of the outcome pairs in order to quantify the three different cases. Nevertheless, almost all students agreed that the largest number of possible pairs corresponded to the most likely event. In particular, they concluded that the grey horse had the highest probability of finishing first after ten rounds, as evidenced by statements extracted from the groups’ A3 sheets (see Fig.~\ref{Cavalli_gruppo1_2}).
\\
It is worth noting that, at this stage, the term \lq\lq probability\rq\rq\ had not yet been formally introduced by the teachers; rather, it was spontaneously used by the students themselves.
\\
These observations are particularly relevant to the research questions guiding the study. On the one hand, they provide insight into the characteristics of mathematical practices that support the introduction of probabilistic concepts, highlighting the role of collective reasoning, comparison of cases, and informal quantification strategies. On the other hand, they inform the identification of design principles for task development, showing how appropriately designed activities can foster the emergence and initial formalization of probabilistic thinking, even in the absence of explicit formal instruction.
With reference to Activity 2, the physical rolling of the two dice showed students that each of the horses could advance during the race; moreover, some races were not won by the grey horse. Nevertheless, the general consensus supported the previously formulated hypothesis that the grey horse had a higher probability of winning, based on the observation that it won the majority of the races conducted.
\\
However, no group proposed a method to quantify the probability of each horse winning based on the observed frequencies. At this stage, one group in particular reported the remark shown in Fig.~\ref{Frequenza_Probabilità}, where a spontaneous connection between the concepts of frequency and probability emerges.
\\
These observations provide further insight into the research questions guiding the study. In relation to the first research question, they highlight the role of experimentation and empirical observation as key characteristics of mathematical practices that support the development of probabilistic reasoning. In relation to the second research question, they suggest that the design of activities involving repeated trials can foster the emergence of connections between frequency and probability, even if students may initially lack the tools to formalize such relationships.
\\
This points to the need for carefully designed instructional interventions that build on these intuitions to support the progressive formalization of probabilistic thinking.
\begin{remark}
It would be interesting, by exploring further \lq\lq situations\rq\rq\ to help students grasp the relationship between the concepts of frequency and probability.
Alternatively, as suggested by the title of the strand \lq\lq Data and Predictions\rq\rq\ one could consider using data (familiar to students) organized into frequency distributions as preparatory work for learning in probabilistic contexts.
\end{remark}
Subsequently, through the observation of the outcomes of the dice pairs, it became easier for almost all groups to organize the counting of the number of pairs corresponding to the three different parity cases (see Fig.~\ref{Coppie_gruppo1_2}).
It was particularly interesting that the students incorporated quantitative assessments into their reasoning as well: they fully calculated the probabilities of the three events using the classical definition (see Fig. ~\ref{Valutazione_probabilità}).
This finding is closely related to the first research question, as it highlights a key characteristic of mathematical practice that can effectively support the introduction of probability in lower secondary education. Furthermore, it addresses the second research question by providing evidence of how carefully designed tasks encourage students to move from informal reasoning to more formalized probabilistic thinking.
\\
For Activity 3, the entire class recognized that the score pairs needed to be distributed so that each horse would have 12 points in its favor. However, no group was able to establish a sufficiently simple rule, as their reasoning continued to focus on individual pairs. At this stage, it became necessary for us to suggest using the sum of the two scores as the basis for the criterion and to remind students how to determine the remainder when an integer is divided by 3.
This observation is closely related to the first research question, as it highlights a key aspect of mathematical practice: students initially rely on concrete, case-by-case reasoning and need guidance to identify more general, systematic strategies that support the introduction of probability concepts. It also addresses the second research question by illustrating how carefully designed instructional interventions can guide students from informal, pair-based reasoning toward more formalized mathematical rules. (See Fig.~\ref{modulo3})
At the conclusion of the first experimentation, several observations can be drawn. In a grade 6th class, students had not yet been introduced to probability concepts. Despite our being unfamiliar to them, the students readily accepted both our presence and their involvement in the proposed activities. Overall, participation was widespread and increased progressively as students became more familiar with the materials and the instructional approach.
At the end of the two-hour session (which extended beyond the scheduled school time), the students asked the following question:
\begin{center}
\lq\lq When can we sit down?\rq\rq
\end{center}
%
%
%
%
\section{Second year laboratory design}
In order to achieve the research objective outlined at the beginning—namely, to derive design principles for the development of probabilistic tasks and learning environments, with particular attention to how these support the gradual formalization of students’ probabilistic reasoning—we returned to the same class in the following year. So that the intervention took place in a second-year lower secondary class at the beginning of February 2021, with the classroom teacher present alongside the students and the researchers participating via videoconference due to the COVID-19 pandemic restrictions in place at the time.
\\
Due to the required physical distancing among students, it was necessary to adopt a collaborative group-work approach in place of cooperative group work. Nevertheless, the use of certain materials previously brought from home by the students (such as dice, boxes, and pre-printed tables) was maintained.
\\
The interactive whiteboard was used to coordinate group work, to better structure the timing of the activities, and to enable interaction from a distance.
\\
The use of materials proved effective in making the game quickly understandable and in highlighting the difference between the outcomes of the throws (pairs or triples) and the resulting sums.
%
%
%
%
%
\subsection{The Game of Zara}
\begin{figure}
    \centering
\includegraphics[width=4.5cm]{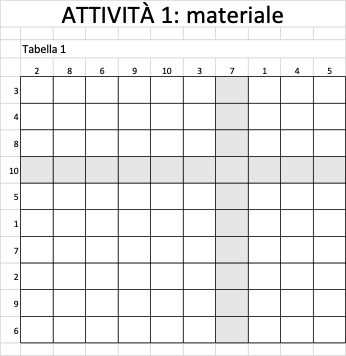}
\caption{Table 1.}
\label{tabella_over}
\end{figure}

{\bf Activity 1}
\begin{itemize}
\item[0)] Place the two dice in the box, cover it with the lid, shake thoroughly, then set the box back on the table. Remove the lid and observe the two numbers on the top faces. Repeat this procedure several times. Examine Table 1 (shown here in Fig.~\ref{tabella_over}) from the materials sheet carefully.
\item[1)] Why are the row labeled 10 and the column labeled 7 shaded in Table 1?
\item[2)] Continue shading the rows in Table 1 that, according to the reasoning identified, are like the row labeled 10.
\item[3)] Continue shading the columns in Table 1 that, according to the reasoning identified, are like the column labeled 7.
\item[4)] How many white cells remain? Explain why.
\item[5)] Rearrange the table on the worksheet so that only white cells remain.
\end{itemize}
\par
\begin{figure}
    \centering
\includegraphics[width=4.0cm]{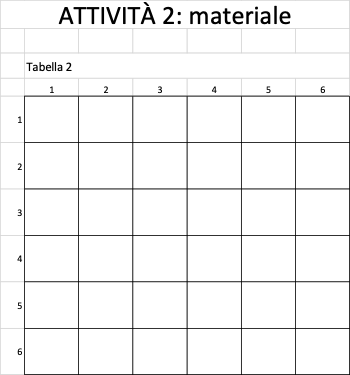}
\includegraphics[width=5.5cm]{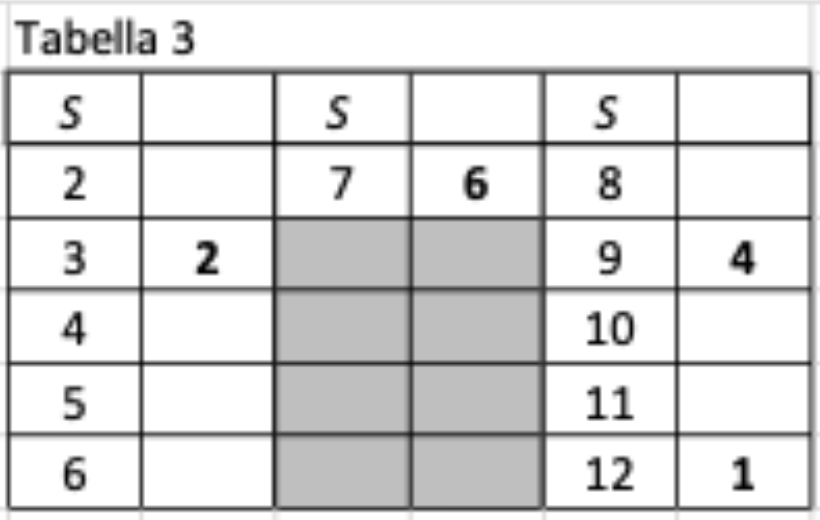}
\caption{On the left, the resulting table showing the relevant data. On the right, the table reporting the frequency of the sums of the outcomes.}
\label{tabelle_finale_somme}
\end{figure}

{\bf Activity 2}
\begin{itemize}
\item[0)] Place the two dice in the box, put on the lid, shake well, and place the box on the table. Remove the lid, observe the two \emph{scores}, and compute their \emph{sum} $S$.
\item[1)] Fill in the cells of Table 2 on the materials sheet (see Fig.~\ref{tabelle_finale_somme}) with the sum $S$ corresponding to the row and column headings.
\item[2)] Identify a reasoning to complete Table 3. Then complete Table 3 (see Fig.~\ref{tabelle_finale_somme}).
\item[3)] With reference to Table 3, what does the expression
$$
(1+2+3+4+5)+6+(5+4+3+2+1)
$$
represent?
\item[4)] How does it differ from the previous expression
$$
2(1+2+3+4+5)+6?
$$
\end{itemize}
\par
{\bf Activity 3}
\begin{itemize}
\item[0)] Place the two dice in the box, put on the lid, shake thoroughly, place the box on the table, then remove the lid. Observe the two scores and compute their sum $S$. Repeat several times. Carefully examine Table 2 in Fig.~\ref{tabelle_finale_somme} on the materials sheet.
\item[1)] We have $6 = 1+5 = 2+4 = 3+3$ and $7 = 1+6 = 2+5 = 3+4$. So, is betting on $S = 6$ as favorable as betting on $S = 7$?
\item[2)] Copy Tables 2 and 3 in Fig.~\ref{tabelle_finale_somme} (completed) into a blank sheet.
\item[3)] Is it more advantageous to bet on $S = 2$ or on $S = 12$? (Justify your answer and indicate which table was more useful.)
\item[4)] Is it more advantageous to bet on $S = 5$ or on $S = 9$? (Justify your answer and indicate which table was more useful.)
\item[5)] Is it more advantageous to bet on $S = 5$ or on $S = 8$? (Justify your answer and indicate which table was more useful.)
\item[6)] On which value of $S$ is it most advantageous to bet?
\end{itemize}
{\bf Activity 4}
\\
This activity was carried out in class by the teacher after a brief \emph{training} session provided by us.
\\
\begin{itemize}
\item[0)] Place the three dice in the box, put on the lid, shake thoroughly, place the box on the table, then remove the lid. Observe the three \emph{scores} and compute their sum $S$. Repeat several times.
\item[1)] What are the possible values of $S$?
\item[2)] It holds: $5 = 1+1+3 = 1+3+1 = 3+1+1 = 1+2+2 = 2+1+2 = 2+2+1$.
\item[3)] It holds: $6 = 1+1+4 = 1+4+1 = 4+1+1 = 1+2+3 = 3+1+2 = 2+1+3 = 2+3+1 = 3+2+1 = 1+3+2 = 2+2+2$.
\item[4)] Is it more advantageous to bet on $S = 5$ or on $S = 6$?
\item[5)] It holds: $16 = 6+6+4 = 6+4+6 = 4+6+6 = 6+5+5 = 5+6+5 = 5+5+6$.
\item[6)] Is it more advantageous to bet on $S = 5$ or on $S = 16$?
\item[7)] In your opinion, on which values of $S$ is it most advantageous to bet?
\end{itemize}
%
%
%
%
%
\subsection{Reflections on the results of the second laboratory}
\begin{figure}
    \centering
\includegraphics[width=10.0cm]{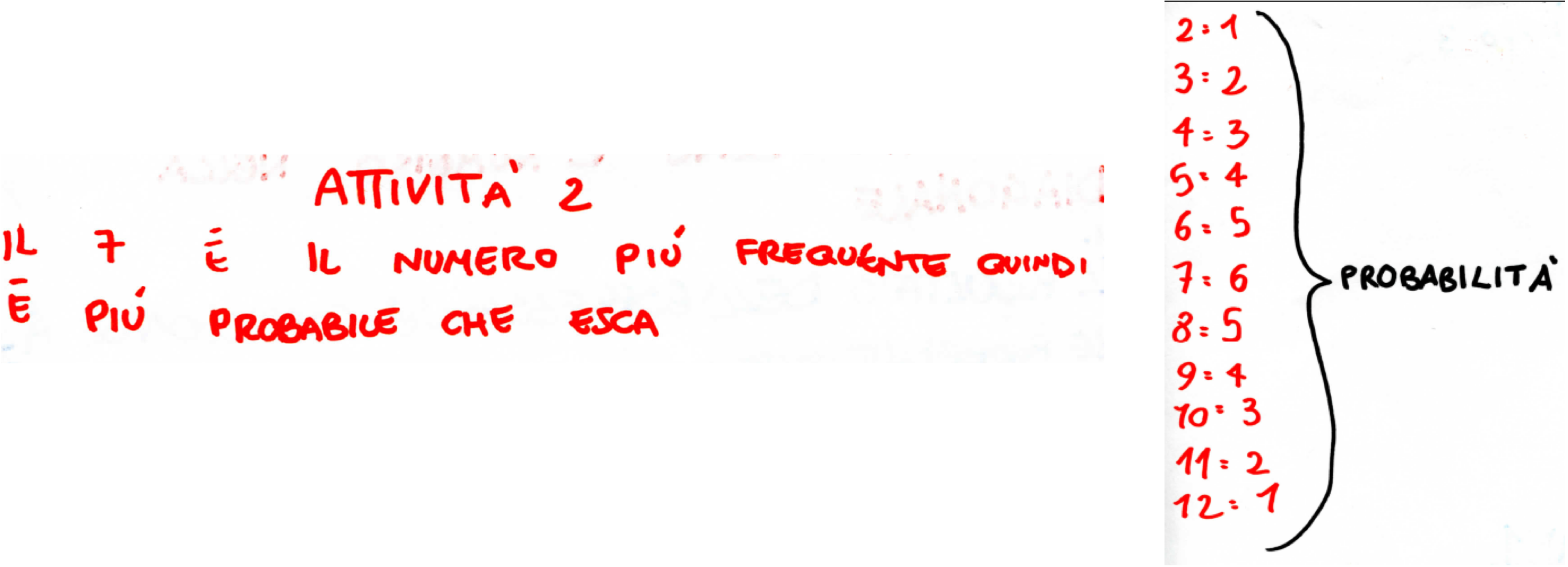}
\caption{On the left \lq\lq 7 is the most frequent number, therefore it is the most probable\rq\rq.}
\label{casi_favorevoli}
\end{figure}

With the preparation provided by the \lq\lq discoveries\rq\rq\ about rolling two dice and the construction of Table 2, in which all 36 pairs are listed, the groups correctly identified the number of favorable cases for each value taken by $S$ (see Fig.~\ref{casi_favorevoli}).
\par
No student asked about the meaning of $S$: they recognized both the meaning and the role of the letter. As a result, all agreed that $S=7$ is the preferable choice for a possible bet (see Fig.~\ref{casi_favorevoli}).
\par
All groups agreed that 7 lies on the secondary diagonal (the \lq\lq longest\rq\rq\ one) of Table 2, but none recognized the role of the diagonals parallel to it for the other values taken by $S$, even though they themselves used different colors to identify the various values of $S$ within Table 2 (see Fig.~\ref{casi_favorevoli}).
\par
\begin{figure}
    \centering
\includegraphics[width=4.0cm]{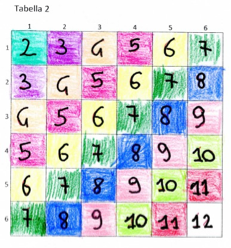}
\caption{Colored Table 2.}
\label{tabella_colorata}
\end{figure}

Only one group gave importance to the central position of 7 among the values taken by $S$.
The idea of \lq\lq preferability\rq\rq\ as a synonym for probability (the ratio between the number of favorable cases and the number of possible cases) was introduced by the students themselves during an oral presentation of their reasoning (Lorenzo).
\\
In particular, with referenze to item 1) in Activity 1, here it is the answer of one student, Giuseppe: \emph{7 and 10 cannot come up on a single die.}. While, with reference to item 4) Lorenzo explains it to his classmates: \emph{here are 36 because we ruled out the others. The teacher suggests looking at the dice, and he starts reasoning like this: if the first die shows a 1, the other die can be any of 1, 2, 3, 4, 5, or 6\ldots\ $6 + 6 + 6 + 6 + 6 + 6$}. 
\\
In the sentence reported in Fig.~\ref{combinazioni} by a group, it is interesting that the word \emph{combinations} is used; although it appears in a counting context, it is not used incorrectly, as it is clear that students are employing it in its everyday sense. In fact, they reinforce its meaning afterward by saying that the numbers \emph{combine}.
\begin{figure}
    \centering
\includegraphics[width=9.0cm]{combinazioni.png}
\caption{Since the highest number that can come up on each die is six, the possible \emph{combinations} are 36. This is because each number can be \emph{combined} with six other numbers.}
\label{combinazioni}
\end{figure}
With reference to items 1) and 2) of Activity 2, we suggested coloring Table 2, as shown in Fig.~\ref{tabella_colorata}, and we asked whether coloring it had been useful.
\\
Federica, who had already written in a column how many times each result appears, said it was not.
\\
Giuseppe and Luigi, on the other hand, said it was useful because it made it easy to count the cells of the same color.
Regarding the last items of Activity 2, it was interesting that, after interpreting the first expression and noticing that the result of the second was the same, even though they used a different property of the operations, they began to write many equivalent expressions by exploiting the properties of the operations.
\par
Concerning the roll of three dice, it was quickly established that there are 216 possible triples and which triples correspond to a given value of $S$. However, most students did not take into account the role of permutations within each triple, reverting to the reasoning known as that of the Florentine gentlemen. This difficulty was reduced following teacher intervention.
\\
No group, however, distinguished the triples according to how many times each die score appeared in them.
\\
Overall, a clear difference was observed between written and verbal communication.
\\
The activity was conducted over two separate days, totaling three hours.
%
%
%
%
%
\section{Conclusions}
The set of competencies considered in this study—ranging from the interpretation of written texts and the use of symbolic representations, to the development of reasoning processes, decision-making based on quantitative information, and the ability to communicate and justify one’s thinking—provides a coherent framework for analyzing and designing mathematical practices in the context of probability education. 
\par
These competencies not only inform the identification of key characteristics of effective mathematical practices that support the introduction of probabilistic concepts in lower secondary school, but also guide the definition of design principles and pedagogical choices underlying the proposed learning activities. 
\par
Through their integration into a sequence of tasks developed and refined over time, they offer a lens for examining students’ engagement with informal probabilistic notions and for tracing the evolution of their reasoning from intuitive interpretations toward progressively more structured and formalized forms. In this perspective, they contribute both to addressing the research questions and to achieving the broader objectives of fostering a meaningful and gradual formalization of probabilistic thinking.
In summary, the integration of design research methodologies with insights from probability education — and with curricular priorities such as those outlined in the Italian national guidelines — provides a powerful foundation for advancing both theory and practice. Such an approach is essential for improving probabilistic literacy and supporting learners in developing the competencies needed to reason effectively under uncertainty.
%
%
%
%
%

%
\end{document}